\documentclass[ECP,preprint]{ejpecp} % replace ECP by EJP if needed.
\usepackage[utf8]{inputenc}

\usepackage[english]{babel}

\SHORTTITLE{Asymptotic behaviour of sampling and transition probabilities in coalescent models }

\TITLE{Asymptotic behaviour of sampling and  transition probabilities in coalescent models under selection and parent dependent mutations } % \thanks is optional. Insert line breaks with \\

\AUTHORS{%
Martina Favero
  \footnote{Stockholm University.
    \EMAIL{martina.favero@math.su.se}}
  \and %% remove this line and below if single author
  Henrik Hult\footnote{KTH, Royal Institute of Technology, Stockholm. \EMAIL{hult@kth.se} }
   }%AUTHORS
 % Separate items with ;

\AMSSUBJ{60J90} % Edit. Separate items with ;
\AMSSUBJSECONDARY{60J70;92D15} % Optional, separate items with ;

\SUBMITTED{ } % Edit.
\ACCEPTED{} % Edit.

\VOLUME{0}
\YEAR{2020}
\PAPERNUM{0}
\DOI{10.1214/YY-TN}

\ABSTRACT{
The results in this paper provide new information on asymptotic properties of classical models: the neutral Kingman coalescent  under a general finite-alleles, parent-dependent mutation mechanism, and its generalisation, the ancestral selection graph.
Several relevant quantities related to these fundamental models are not explicitly known when mutations are parent dependent. Examples include the probability that a sample taken from a population has a certain type configuration,  and the transition probabilities of their block counting  jump chains.
In this paper,  asymptotic results are derived for these quantities,
as the sample size goes to infinity.
It is shown that the sampling probabilities decay polynomially in the sample size with multiplying  constant depending on the stationary density of the Wright-Fisher diffusion   
and that the transition probabilities  converge to the limit of frequencies of types in the sample.
}

\newcommand{\E}[1]{\mathbb{E}\left[#1\right]}
\newcommand{\Prob}[1]{\mathbb{P} \left( #1\right) }
\newcommand{\x}{\textbf{x}}
\newcommand{\y}{\textbf{y}}
\newcommand{\n}{\textbf{n}}
\newcommand{\e}{\textbf{e}}
\newcommand{\X}{\textbf{X}}
\newcommand{\Y}{\textbf{Y}}
\newcommand{\norm}[1]{\left\lVert#1\right\rVert_1}

\newcommand{\nth}{^{(n)}}
\newcommand{\domain}{\Delta}
\newcommand{\kc}{\mathbf{H}}
\newtheorem{Remark}[theorem]{Remark}
\begin{document}

%%%%%%%%%%%%%%%%%%%%%%%%%%%%%%%%%%%%%%%%%%%%%%%%%%%%%%%%%%%%%%%%%%%
%%                                                               %%
%% No need for \maketitle.                                       %%
%%                                                               %%
%%%%%%%%%%%%%%%%%%%%%%%%%%%%%%%%%%%%%%%%%%%%%%%%%%%%%%%%%%%%%%%%%%%

%%%%%%%%%%%%%%%%%%%%%%%%%%%%%%%%%%%%%%%%%%%%%%%%%%%%%%%%%%%%%%%%%%%
%%                                                               %%
%% Please replace what follows by the body of your article       %%
%% (up to the bibliography):                                     %%
%%                                                               %%
%%%%%%%%%%%%%%%%%%%%%%%%%%%%%%%%%%%%%%%%%%%%%%%%%%%%%%%%%%%%%%%%%%%

\section{Introduction}
\label{sect:intro}
Given a sample of genetic material from some individuals in a neutral population, the Kingman coalescent \cite{kingman1982b} models the genealogy of the individuals.
Under the finite-alleles model of mutation, when mutations are parent independent (PIM), i.e. the genetic type of the mutated offspring does not depend on the type of its parent, the coalescent can be easily simulated backwards in time, i.e. starting from the individuals in the sample and going back to their most recent common ancestor, since its backward transition probabilities are explicitly known. Furthermore,  the likelihood function, also referred to as sampling probability, is explicitly known as well as several other quantities. However, when mutations are parent dependent, the model becomes more complex and several of the expressions, that are explicit in the PIM case, become implicit. 

The same dichotomy occurs when, in addition to mutation, selection is considered and the Kingman coalescent is replaced by its well-known generalisation: the ancestral selection graph (ASG) \cite{Krone1997,Neuhauser1997}.
The ASG represents the history of a sample not only through mutation and coalescence events, but also through branching events. A branching of a lineage, into a true and a virtual lineage, corresponds to a selection event: two individuals are chosen as potential parents, the most viable of them becomes the true parent.

Coalescent models, i.e. the Kingman coalescent and its numerous generalisations,  have been extensively used for inference on genetic data sets, in combination with Monte Carlo methods used to approximate implicit likelihood functions in the non-PIM case,  %e.g. \cite{stephens2000,deiorio2004,Koskela2015,Koskela2018}, 
see e.g. \cite{Stephens2007} for an overview. 
The asymptotic behaviour of coalescent models in relation to samples of large size has recently gained attention due to the large size of modern study samples. 
While the size of samples continues to increase, due to advancements in DNA sequencing technology, inference methods based on coalescent models struggle to provide reliable results, since the coalescent does not scale well in terms of sample size \cite{kelleher2016}. Moreover, as discussed in \cite{bhaskar2014}, there is a distortion of some of the properties of the coalescent, which is a suitable approximation of certain models, e.g. Wright-Fisher model, provided the sample size is sufficiently smaller than the effective population size. In some cases this leads to inaccurate conclusions, e.g. concerning prediction of rare variants. 
Therefore, 
studying large-sample-size properties of coalescent models is not only an interesting theoretical problem, but also provides tools for addressing the above mentioned shortcomings, which are directly related to practical applications. 
More precisely,
because of the large size of modern samples, inference methods  based on coalescent models are widely used, 
even when the underlying assumption of sample size being sufficiently smaller than the effective sample size is violated, as discussed in \cite{bhaskar2014}. Thus a theoretical large-sample-size efficiency analysis of algorithms based on coalescent models, such as the ones in \cite{stephens2000,stephens2003,deiorio2004,birkner2008,griffiths2008,hobolth2008,birkner2011,Koskela2015,Koskela2018}, would be useful. The study of large-sample-size asymptotic  properties  of the coalescent, to which this paper and \cite{favero2020c} aim to  contribute, are relevant for such analysis. 
Furthermore, large-sample-size  asymptotic results provide tools for the  analysis of differences between the coalescent  approximation and the original model. In fact,  a direct theoretical  comparison of the properties of these models is challenging, in \cite{bhaskar2014} a numerical comparison is provided,  whereas comparing the corresponding simpler  limiting objects may provide interesting insights.   
%could likely serve also as a support in the analysis of these problems.

This paper provides an analysis of the asymptotic behaviour of 
the sampling probabilities under the ASG, and of the transition probabilities of its block counting jump chain. Note that the Kingman coalescent can be seen as a special case of ASG by setting the selection parameters equal to zero.
A finite-alleles model for mutation is assumed, with $d$ possible types.

A sample is represented by a vector $\n\in \mathbb{N}^d\setminus\{0\}$, where $n_i$ represents the number of individuals in the sample carrying allele $i$, $i=1,\dots,d$.
The likelihood function, or sampling 
probability, $p(\n)$, corresponds to the probability that a sample taken from a population at equilibrium has a  configuration of types given by the vector $\n$.  
A recursion formula is known for $p$, see  \cite{Krone1997} for the ASG,  \cite{griffiths1994simulating} and the references therein for the particular case of the Kingman coalescent. However, an explicit formula is unknown in the general case of parent dependent mutations. When the sample size, $\norm{\n}=n_1+\dots + n_d$, is large,  it is computationally too expensive to compute $p$ using the recursion formula.
In this paper a large sample of the form 
$n \y\nth$, with $n\in\mathbb{N}\setminus\{0\}$ and 
$\y\nth\in \frac{1}{n}\mathbb{N}^d\setminus\{\boldsymbol{0}\}$, is considered and
the asymptotic behaviour of the sampling probabilities $p(n\y\nth)$ is studied,
as $n$ tends to infinity and $\y\nth$ tends to some $\y\in \mathbb{R}_{>0}^d$.
In particular, 
 it is proved that the sampling probabilities decay polynomially, that is,
 \begin{equation}
    \label{eq:asymptp}
    p(n\y\nth)\sim  \tilde{p}\left(\frac{\y}{\norm{\y}}\right)\norm{\y}^{1-d}
    n^{1-d}, \quad \text{as } n \to \infty, \ \y\nth\to\y,  
    \end{equation}
where $\sim$ denotes asymptotic equality in the sense that $a_n \sim b_n$ if $\lim_{n \to \infty} a_n/b_n = 1$, $\norm{\y} = |y_1| + \dots + |y_d|$, and $\tilde{p}$ is the stationary density of the Wright-Fisher diffusion that is dual to the ASG, see Section \ref{sect:framework} for more details on $\tilde{p}$. 

The intuition for the polynomial decay comes from the neutral scenario, when mutations are parent independent, in fact, in this case the sampling probabilities are explicitly known, see e.g. \cite{griffiths1994simulating}, and it is straightforward to show that their asymptotic decay is polynomial of degree $d-1$, 
see Subsection \ref{subsect:PIM} for an explicit calculation.
In particular, the degree of the polynomial does not depend on the mutation parameters, although the multiplicative constant does, which hints that the same behaviour applies in general, even when mutations are parent dependent.
To address the general case the following strategy is adopted.  
The classical representation formula of the sampling probabilities as  
expectations with respect to the stationary distribution of the Wright-Fisher diffusion, which is recalled in Section \ref{sect:framework},  is  used in Subsection \ref{subsect:interpret} to  interpret the sampling probabilities as  expectations with respect to a sequence of Dirichlet distributions. A local limit theorem for the sequence of Dirichlet distributions is derived in the Appendix. Finally, in Section \ref{sect:pconv} the general result is proved.
   
Establishing the asymptotic decay of $p(n\y\nth)$, enables the study of the asymptotic behaviour of the transition probabilities of the block counting jump chain of the 'typed' ancestral selection graph, to which Section \ref{sect:probabilities} is dedicated.
The adjective 'typed' is used to put emphasis on the fact that, in this paper, the lineages of the ASG  are always associated to a type, as in e.g. \cite{Etheridge2009,favero2021}, 
whereas often, e.g. in \cite{Krone1997}, the ASG is first constructed backwards in time with no regard for types, which are superimposed afterwards on the graph.
This two-steps construction does not allow to construct the genealogy of a sample with given types, which is of interest when performing inference based on a sample. 
The block counting process of the 'typed' ASG is the process that counts how many lineages of each type are present as time evolves from when the sample is taken until the most recent common ancestor is reached. In this paper the jump chain of this block counting process is considered, and the asymptotic behaviour of its transition probabilities is studied in Section \ref{sect:probabilities}, where the chain is properly scaled.

%%%%%%%%%%%%%%%%%%%%%%%%%%%%%%%%%%%%%%%%%%%%%%%%%%%%%%%%%%%%%%%%%%%%%%%%%%%%%%%%%%
\section{Framework}
\label{sect:framework}

In Section \ref{sect:intro}, $p(\n)$ is defined as the probability of sampling $\n$ from a population at equilibrium, when the underlying model is the ASG. In this section, the ASG, its block counting jump chain and the related Wright-Fisher diffusion are introduced; furthermore, the properties that are relevant for the results in this paper are recalled.

To set the notation and recall the definition of the model, 
let $\theta$ be the mutation rate, $P=(P_{ij})_{i,j=1}^d$ be the mutation probability matrix and $\gamma_i,i=1,\dots,d$, be negative selection parameters. Note that the selection parameters can also be chosen to be equal to zero, in that case the ASG becomes the classical Kingman coalescent.

The ASG describes the evolution of lineages over time. When $m$ lineages, i.e. edges, are present, either one of them is lost by coalescence with another, 
at rate $\binom{m}{2}$, or one is added, by branching of an existing lineage, at rate $m\frac{\gamma }{2}$, where $\gamma=\max\{\gamma_i-\gamma_j:i,j=1,\dots,d\}$. 
Furthermore, 
one of the $m$ lineages undergoes a mutation event at rate $m\frac{\theta}{2}$. 
The evolution continues, for an almost surely finite time, until only one lineage is left, representing the most common ancestor of the initial $m$ lineages.
To assign types to each lineage, a type is assigned to the ancestor and accordingly to all the descending lineages by following the graph structure, which includes mutation points, and using the mutation matrix, $P$, and the selection parameters, $\gamma_1,\dots,\gamma_d$, to determine the type of a lineage after mutation and selection (branching) events.
For a complete and rigorous definition of the ASG, see \cite{Krone1997,Neuhauser1997}.

In this paper it is enough to consider only the jump chain of the block counting process of the ASG, instead of the entire graph structure, and it is crucial to assign types to lineages while the 'typed' chain evolves backwards in time, as opposed to the a-posteriori type assignment described above. This enables the construction of the genealogy of a 'typed' sample, at the cost of an implicit backward transition distribution.
The block counting process of the 'typed' ASG is described in detail in \cite{Etheridge2009}. Its jump chain, $\kc=\{\kc(l)\}_{l\in\mathbb{N}}\in \mathbb{N}^d\setminus\{\boldsymbol{0}\}$, which is the focus of Section \ref{sect:probabilities},  
is the time homogeneous Markov chain with the following transition probabilities \cite{Etheridge2009}, for $\n\in\mathbb{N}^d\setminus\{\boldsymbol{0}\}$, 
    \begin{equation}
    \label{defrho}
    \begin{aligned}
    %=&
    &\Prob{\mathbf{H}(l+1)=\n-\mathbf{v}
    \mid
    	\mathbf{H}(l)=\n}
    %\\
    =p(\n-\mathbf{v}\mid\n)
    =\\
    &=\begin{cases}
    \frac{n_j( n_j -1)}{\sum_{r=1}^d n_r|\gamma_r|+\norm{\n}( \norm{\n} -1 + \theta )} \frac{1}{\pi[j|\n-\e_j]}, \!\!\!\!
    &\text{  if  }
    \mathbf{v}=\mathbf{e}_j,
    \; j\! =\! 1\dots d,
    \\
    \frac{\theta P_{ij} n_j }
    {\sum_{r=1}^d n_r|\gamma_r|+ \norm{\n} ( \norm{\n}-1 + \theta ) }
    \frac{\pi[i|\n-\e_j]}{\pi[j|\n-\e_j]},\!\!\!\!
    &\text{  if  }
    \mathbf{v}=\mathbf{e}_j-\mathbf{e}_i,
    \; i, j\! = \! 1\dots d,
    \\
    \frac{\norm{\n} |\gamma_j| }
    {\sum_{r=1}^d n_r|\gamma_r|+ \norm{\n} ( \norm{\n}-1 + \theta ) }
    \pi[j|\n],\!\!\!\!
    &\text{  if  }
    \mathbf{v}=-\mathbf{e}_j,
    \;  j\! = \! 1\dots d,
    \\
    0, \!\!\!\! &\text{  otherwise,}
    \end{cases}
    \end{aligned}
    \end{equation}
where $\e_j$ is the $j^{th}$ $d-$dimensional unit vector and $\pi[j|\n]$ is the  probability of sampling an individual of type $j$ after sampling $\norm{\n} $ individuals with types given by $\n$.
Note that a step $+\e_j$, resp. $-\e_j$,  corresponds to the coalescence, resp. branching, of a lineage of type $j$, and $\e_j-\e_i$ corresponds to the mutation of a lineage from type $i$ to type $j$.
The probability $\pi$
can be written equivalently in terms of the sampling probabilities, $p$, as 
    \begin{align}
    \label{defpi}
    \pi[i|\n]= \frac{n_i+1}{\norm{\n}+1} \frac{p(\n+\e_i)}{p(\n)} .
    %=\frac{k(\n+\e_i)}{k(\n)}
    \end{align}
%see \cite{stephens2000,deiorio2004} for more details.
Similarly to the sampling probabilities, the probability $\pi$ is unknown explicitly, unless mutations are parent independent.

The transition probabilities, which are the focus of Section \ref{sect:probabilities}, are related to the recursion formula for $p(\n)$ mentioned in the introduction. 
In some cases, as in this paper, instead of using the  recursion formula, it is convenient to work with a well-known representation for $p(\n)$ in terms of the Wright-Fisher diffusion, which is described in the following.  

While the ASG models the ancestral history of a sample taken from the population, the allele frequencies are modelled by the Wright-Fisher diffusion 
    $
    \X=\{\X(t)\}_{t\geq 0}\subset 
    \mathcal{S}=\{\x\in[0,1]^{d}:\sum_{i=1}^{d}x_i=1\}
    $
that is the solution to the following stochastic differential equation, 
    \begin{align}
    \label{eq:sdeWF}
    d\X(t)= \mu(\X(t)) dt + \sigma(\X(t))^{1/2} d \mathbf{W}(t), \quad t \geq 0, 
    \end{align}
where $\mathbf{W}=\{\mathbf{W}(t)\}_{t\geq 0}$ is a $d$-dimensional Wiener process,
the diffusion matrix $\quad$ is 
    $
    \sigma_{ij}(\x)=
    x_i(\delta_{ij}-x_j),
    i,j=1,\dots,d
    $,
and the drift is
    $
    \mu_i(\x)=
    \theta \sum_{j=1}^{d} x_jP_{ji}- \theta x_i
    +x_i\left( \gamma_i- \sum_{j=1}^d\gamma_j x_j\right)
    ,i=1,\dots,d
    $,
see e.g. \cite{Etheridge2011}. When the mutation probability matrix $P$ is irreducible, the stationary distribution of the Wright-Fisher diffusion $\X$ is unique and has a smooth density, $\tilde{p}$, with respect to the Lebesgue measure, to which Remark \ref{remark:stationarity} is dedicated.
In this paper we assume that $P$ is irreducible.  

The sampling probability, $p(\n)$, can be expressed as the expectation of a multinomial draw from the stationary distribution of the Wright-Fisher diffusion. Let $\tilde{\X}$ be distributed according to  the stationary distribution of the Wright-Fisher diffusion, then  
    \begin{align}
    \label{eq:multinomialdraw}
    p(\n) = 
    \binom{\norm{\n}}{\n}
    \E{\prod_{i=1}^{d}\tilde{X}_i^{n_i}}.
    \end{align}
This formula is a consequence of the duality relationship between the ASG and the diffusion, in fact, when such relationship is proved, it is also shown that the right hand side of \eqref{eq:multinomialdraw} solves the recursion formula that defines the sampling probability, see for example  \cite{Krone1997,favero2021}.
Furthermore,
since the sample is exchangeable, the formula can also be explained by de Finetti's representation theorem.

\begin{Remark}[Stationary density of Wright-Fisher diffusion]
\label{remark:stationarity}
Equation \eqref{eq:multinomialdraw} only holds if the Wright-Fisher diffusion admits a stationary distribution. 
The existence of a unique stationary distribution is related to the structure of the mutation mechanism, more precisely,  to the existence of an invariant measure for the mutation probability matrix $P$. 
When mutations are parent independent, $P_{ij}=Q_j>0$, $i,j=1,\dots,d$, 
the stationary distribution, not only exists, but also has an explicitly known density:
in the neutral case 
a Dirichlet density with parameters $\theta Q=(\theta Q_1,\dots,\theta Q_d)$, see e.g.  \cite{Wright1949,griffiths1994simulating,Etheridge2009}, and when selection is included, 
 a weighted Dirichlet density, see e.g.  \cite{Etheridge2009,favero2021}.
Unfortunately, the PIM case is the only case where the stationary distribution is explicitly known.
Furthermore, the Wright-Fisher diffusion $\X$ has a degenerate elliptic generator, since the diffusion matrix has zero entries when one of the components of $\X$ is equal to zero. 
%on the boundary of the domain. 
Thus the classical theory for the study of stationarity and  for the study of solutions to the Fokker-Planck equation does not apply and an ad-hoc analysis, which is provided in \cite{shiga1981}, is needed. 
In \cite[Thm 3.1]{shiga1981} it is shown that the stationary distribution of $\X$ exists uniquely, 
assuming that  $P$ is irreducible, which  implies the invariant measure of $P$ exists uniquely.
Furthermore, in \cite[Thm 3.2]{shiga1981} it is also shown that, under the same assumption, the stationary distribution is absolutely continuous with respect to the Lebesgue measure and its probability density function, $\tilde{p}$, defined on $\domain=\{\x\in[0,1]^{d-1}:\sum_{i=1}^{d-1}x_i\leq 1\}$, 
is smooth on 
$\domain^{\mathrm{o}}=\{\x\in(0,1)^{d-1}:\sum_{i=1}^{d-1}x_i< 1\}$.
\end{Remark}

In the light of the previous remark, in this paper it is assumed that $P$ is irreducible,
in fact, the existence of a smooth stationary density, even of an unknown form, proves to be sufficient for studying the asymptotic behaviour of the sampling probabilities through \eqref{eq:multinomialdraw}.  

%%%%%%%%%%%%%%%%%%%%%%%%%%%%%%%%%%%%%%%%%%%%%%%%%%%%%%%%%%%%%%%%%%%%%%%%%%%%%%%%%

\section{Sampling probabilities}
\label{sect:isamplingprob}

The aim of this section is to provide a representation of the sampling probabilities that is convenient for the study of  their asymptotic behaviour. This is a key step, in fact, once the representation is identified, 
the outline of the proof of the asymptotic result  becomes intuitively clear.

\subsection{Parent independent mutations under neutrality}
\label{subsect:PIM}

Before focusing on the general case,  the  PIM  case without selection is analysed in order to provide better insight. Assume $P_{ij}=Q_j, \gamma_j=0, i,j=1,\dots,d$ and let $Q=(Q_1,\dots,Q_d).$ 
As explained in Remark \ref{remark:stationarity}, in this case the stationary density is Dirichlet with parameters $\theta Q$, thus,  
computing the expectation in \eqref{eq:multinomialdraw}, 
%and the Dirichlet stationary density for the PIM neutral case, 
yields the following explicit expression for  the sampling probabilities
    $$
    p(\n)= \binom{\norm{\n}}{\n} \frac{B(\n+\theta Q )}{B(\theta Q)}
    =
    \frac{1}{B(\theta Q)}
    \frac{\Gamma(\norm{\n}+1 )}{\Gamma(\norm{\n}+\theta)}
    \prod_{i=1}^{d}\frac{\Gamma(n_i+\theta Q_i )}{\Gamma(n_i+1)},
    $$
where $B$ is the multidimensional Beta function and $\Gamma$ is the Gamma function, the expression above can be found in e.g. \cite{griffiths1994simulating}. Applying  Stirling's formula to the Gamma functions yields, as $\y\nth\to\y,$
    \begingroup
    \allowdisplaybreaks
    \begin{align*}
    p(n\y\nth)
    &\sim
    \frac{1}{B(\theta Q)}
    \frac{
	\left(\norm{n\y\nth}+1 \right)^{\norm{n\y\nth}+\frac{1}{2}} e^{-\left(\norm{n\y\nth}+1\right)}
    }{
    \left(\norm{n\y\nth}+\theta\right)^{\norm{n\y\nth}+\theta-\frac{1}{2}}e^{-\left(\norm{n\y\nth}+\theta\right)}
    }
    \\
    &\quad\quad\quad
    \cdot\prod_{i=1}^{d}
    \frac{
	\left(n y\nth_i+\theta Q_i \right)^{n y\nth_i+\theta Q_i-\frac{1}{2}}e^{-\left(n y\nth_i+\theta Q_i\right)}
    }{
    \left(n y\nth_i+1\right)^{n y\nth_i+\frac{1}{2}}e^{-\left(n y\nth_i+1\right)}
    }
    \\
    &\sim
    \frac{1}{B(\theta Q)}
    \left(1+ \frac{1-\theta}{\norm{n\y\nth}+\theta} \right)^{\norm{n\y\nth}+\frac{1}{2}} 
    \left(\norm{n\y\nth}+\theta\right)^{1-\theta} e^{\theta -1}
    \\
    &\quad\quad\quad
    \cdot\prod_{i=1}^{d}
    \left(1+ \frac{\theta Q_i -1}{n y_i\nth+\theta Q_i} \right)^{n y_i\nth+\frac{1}{2}} 
    \left(n y_i\nth+\theta Q_i\right)^{\theta Q_i -1} e^{1- \theta Q_i}
    \\
    &\sim
    n^{1-d}
    \norm{\y}^{1-d}
    \frac{1}{B(\theta Q)}
    \prod_{i=1}^{d}\left(\frac{y_i}{\norm{\y}}\right)^{\theta Q_i-1}
    \end{align*}
    \endgroup
Note that the sampling probabilities decay polynomially with degree $d-1$. While the multiplicative constant depends on the mutation parameters, the degree,  $d-1$, does not. For this reason one may expect the same degree of decay when mutations are parent dependent, 
which is indeed  correct, and not affected by selection, as shown in Section \ref{sect:pconv}. Note also that the limiting behaviour in the last display corresponds to what was anticipated in \eqref{eq:asymptp}, since the stationary density in this case is a Dirichlet density. 

\subsection{Interpreting the sampling probabilities}
\label{subsect:interpret}
Despite the ease of the intuition, the proof of an asymptotic result for $p(n\y\nth)$ in the general case is more involved, because of the lack of an explicit form for the stationary density of the Wight-Fisher diffusion.
To study the asymptotic behaviour of the normalised probabilities $n ^{d-1} p(n\y\nth)$, it is tempting to interchange the limit and integration in \eqref{eq:multinomialdraw} and consider
    \begin{align*}
     \int_{\domain} \lim_{n \to \infty} n ^{d-1}
    \binom{\| n\y\nth \|_1}{n\y\nth} \prod_{i=1}^d x_i^{n y\nth_i} \tilde{p}(\x) d\x,  
    \end{align*}
The difficulty with this approach is to justify interchanging the limit and integration. Moreover, as $n\to\infty, \y\nth\to\y$, one can show, using Stirling's approximation, that the integrand approaches a Dirac delta function at $\y$. Consequently, the limit is not well defined and care must be taken to rigorously prove the asymptotic behaviour. To this end, an alternative representation of the sampling probabilities is preferred, as outlined below. 

Denote 	the expectation in \eqref{eq:multinomialdraw} by $k(\n)$, i.e., for $\n \in \mathbb{N}^d\setminus\{\boldsymbol{0}\}$,
    \begin{align}
    \label{eq:kdef}
    k(\n)=\E{\prod_{i=1}^{d}\tilde{X}_i^{n_i}}
    =\int_{\domain} 
    \prod_{i=1}^{d} x_i^{n_i} \tilde{p}(\x) d\x, 
    \end{align}
so that, by \eqref{eq:multinomialdraw}, for  $\y\nth \in \frac{1}{n}\mathbb{N}^d\setminus\{\boldsymbol{0}\}$,
    \begin{equation}
    \label{eq:p_intermsof_k}
    p(n\y\nth) = \binom{ \norm{n\y\nth}}{ n\y\nth} k(n\y\nth).  
    \end{equation}
The function $k$ admits an interpretation as an expectation with respect to a Dirichlet distribution, provided it is divided by the appropriate normalising constant, i.e.
$B(n\y\nth+\boldsymbol{1})$, where $\boldsymbol{1}$ is the vector of ones in $\mathbb{R}^d$. By identifying 
    $
    f_{\mathbf{D}\nth}(\x):=B(n\y\nth+\boldsymbol{1})^{-1}
    \prod_{i=1}^{d}x_i^{n_i}
    $
as the probability density function of a Dirichlet random vector, $\mathbf{D}\nth$, with concentration parameters equal to $n\y\nth+\boldsymbol{1}$, \eqref{eq:kdef} is equivalent to the following expression
    \begin{align}
    \label{eq:kexpdirich}
    \frac{ k(n\y\nth)}{ B(n\y\nth+\boldsymbol{1})}
    =
    \int_{\domain} 
    f_{\mathbf{D}\nth}(\x) \tilde{p}(\x)d\x =
    \E{\tilde{p}(\mathbf{D}\nth)}.
    \end{align}

A clarification about notation and state spaces is needed. The space $\mathcal{S}$ is used as a state space for the Wright-Fisher diffusion, $\X$, and the Dirichlet distributed random vectors, $\mathbf{D}\nth$, while their densities are integrated over $\domain$.  
When the density functions, which are defined on $\domain$, are evaluated at a point in $\mathcal{S}$, it is implicitly understood that the first $d-1$ components of the vectors are used, the last component being a function of the first $d-1$ components. And vice versa, for a vector in $\x\in\domain$, $x_d$  stands for $1-\sum_{i=1}^{d-1}x_i$. 

Interpreting the function $k$ as in \eqref{eq:kexpdirich}, 
turns out to be an effective tool for 
analysing the asymptotic behaviour of the sampling probabilities. 
In fact, it can be proved that the expectation in \eqref{eq:kexpdirich}  converges to the constant $\tilde{p}\left(\frac{\y}{\norm{\y}}\right)$ (see Theorem \ref{thm:kconv} for details),
while the remaining factors, 
    $ 
    \binom{n\norm{\y\nth}}{n\y\nth}
    $ 
    $
    B(n\y\nth+\boldsymbol{1})
    $, 
give rise to the polynomial decay of $p(n\y\nth)$ (see Theorem \ref{thm:asymptoticp} for details). 

%%%%%%%%%%%%%%%%%%%%%%%%%%%%%%%%%%%%%%%%%%%%%%%%%%%%%%%%%%%%%%%%%%%%%%%%%%%%%%%%%%%%

\section{Asymptotic behaviour of the sampling probabilities}
\label{sect:pconv}

The study of the asymptotic behaviour of the Dirichlet random vectors $\mathbf{D}\nth$ 
appearing in \eqref{eq:kexpdirich} is reported in the Appendix and is summarized in the following proposition.

\begin{proposition}
\label{prop:convdirich}
Let 
    $
    \mathbf{D}\nth\sim \text{Dirich}(\boldsymbol{\alpha}\nth)
    $, 
	$
	\boldsymbol{\alpha}\nth\in \mathbb{R}_{\geq0}^d
	$ 
such that
	$
	\lim_{n\to\infty}\frac{\boldsymbol{\alpha}\nth}{n}=\boldsymbol{\alpha}\in \mathbb{R}_{>0}^d.
	$
Then the following central limit theorem holds
	$$
	\sqrt{n}\left(
	\mathbf{D}\nth 
	-  \frac{\boldsymbol{\alpha}}{\norm{\boldsymbol{\alpha}}}
	\right) \xrightarrow[n\to\infty]{d}
	\mathcal{N}_d
	\left(\boldsymbol{0},\Sigma(\boldsymbol{\alpha})\right),
	$$
with
	$$
	\Sigma_{ij}(\boldsymbol{\alpha})= 
	\frac{\alpha_i}{\norm{\boldsymbol{\alpha}}^3}
	(\delta_{ij}\norm{\boldsymbol{\alpha}}-\alpha_j),
	\quad i,j=1,\dots,d .
	$$
Furthermore, a local limit theorem for the corresponding probability  density functions, 
$\phi_n$ and $\phi$, holds 
	$$
	\lim_{n\to\infty}\sup_{\mathbf{u}\in\mathbb{R}^{d-1}}|\phi_n(\mathbf{u})-\phi(\mathbf{u})|=0.
	$$
\end{proposition}
In order to prove  an asymptotic result for the sampling probabilities,  the asymptotic behaviour of the expectations \eqref{eq:kexpdirich} is studied in the next theorem, 
using Proposition \ref{prop:convdirich}, 
the  main difficulty  being that the stationary density is unknown and possibly unbounded near the boundary.
\begin{theorem}
\label{thm:kconv}
Let $k$ be defined as in \eqref{eq:kdef}.
Let $\tilde{p}$ be the stationary density of the Wright-Fisher diffusion 
\eqref{eq:sdeWF}, assuming the mutation probability matrix, $P$, is irreducible. 
Let $\y\nth\in \frac{1}{n}\mathbb{N}^d\setminus\{\boldsymbol{0}\}$ such that $\y\nth\to \y\in\mathbb{R}^d_{>0}$, as $n\to\infty$.
Then
	$$
	\lim_{n\to\infty} 
	\frac{k\left(n\y\nth\right)}{B\left(n\y\nth+\boldsymbol{1}\right)}
	=
	\tilde{p}\left(\frac{\y}{\norm{\y}}\right).
	$$
	
\begin{proof}
As explained in Remark \ref{remark:stationarity}, by \cite{shiga1981}, the stationary density $ \tilde{p}$ is smooth on  $\domain^{\mathrm{o}}$, and thus  bounded on any compact set contained in $\domain^{\mathrm{o}}$. It could, however, explode on the boundary. In order to deal with this problem, the domain is divided in two parts. 
For $\varepsilon>0,$
define
	$
	\domain^\varepsilon=
	\{\x\in\domain^{\mathrm{o}}:x_i\geq \varepsilon,i=1,\dots,d\}
	$.
Since
	$
	\frac{\y}{\norm{\y}}\in \domain^{\mathrm{o}}
	$,
it follows that 
	$
	\frac{\y}{\norm{\y}}\in \domain^\varepsilon
	$,
for all $0<\varepsilon\leq \varepsilon_\y$,  
with $\varepsilon_\y=\frac{1}{\norm{\y}}\min_{i=1,\dots,d}y_i$.
Fixing $0<\varepsilon\leq\varepsilon_\y $, rewrite \eqref{eq:kexpdirich} as
	\begin{align}
	\label{eq:intsplit}
	\frac{k\left(n\y\nth\right)}{B\left(n\y\nth+\boldsymbol{1}\right)}
	&=
	\int_{\domain} f_{D\nth}(\x)\tilde{p}(\x)d\x \nonumber 
	\\& =
	\int_{\domain^\varepsilon} f_{D\nth}(\x)\tilde{p}(\x)d\x
	+
	\int_{\domain\setminus \domain^\varepsilon} f_{D\nth}(\x)\tilde{p}(\x)d\x
	.
	\end{align}
To show convergence for the first term in the RHS of \eqref{eq:intsplit},  a change of variables yields
	$$
	\int_{\domain^\varepsilon} f_{D\nth}(\x)\tilde{p}(\x)d\x
	=
	\int
	\mathbb{I}_{\sqrt{n}\left(\domain^\varepsilon- \frac{\y\nth}{\norm{\y\nth}}\right)} (\mathbf{u})
	\phi_n (\mathbf{u}) \tilde{p}\left(\frac{1}{\sqrt{n}}\mathbf{u}+\frac{\y\nth}{\norm{\y\nth}}\right)
	d\mathbf{u},
	$$
where $\phi_n$ and $\phi$ are the density functions defined in Proposition \ref{prop:convdirich} with   $\boldsymbol{\alpha}\nth=n\y\nth+\boldsymbol{1}$. 
By Proposition \ref{prop:convdirich} and continuity of $\tilde{p} $ on $\domain^{\mathrm{o}}$, 
the integrand above converges pointwise to 
	$
	\phi (\mathbf{u})
	\tilde{p}\left(\frac{\y}{\norm{\y}}\right).
	$
Furthermore, 
since $\tilde{p}$ is bounded in $\domain^\varepsilon$ by some $c_\varepsilon$,
the sequence is dominated by the sequence $\phi_n(\mathbf{u})c_\varepsilon$, the integral of which is equal to $c_\varepsilon$, for all $ n$. 
Therefore, by the general dominated convergence theorem 
\cite{royden2010}, 
	$$
	\int_{\domain^\varepsilon} f_{D\nth}(\x)\tilde{p}(\x)d\x
	\xrightarrow[n\to\infty]{}
	\int_{\mathbb{R}^{d-1}}
	\phi(\mathbf{u})\tilde{p}\left(\frac{\y}{\norm{\y}}\right)d\mathbf{u}
	= \tilde{p}\left(\frac{\y}{\norm{\y}}\right).
	$$
It remains to show that 
the second term in the RHS of \eqref{eq:intsplit} converges to zero.
First note that, if $\x\in\domain\setminus\domain^\varepsilon$, then $x_j<\varepsilon$ for some $j$, and, since $x_i\leq 1,i=1,\dots,d$,
	$$
	\prod_{i=1}^{d}x_i^{ny_i\nth}
	< \varepsilon^{ny_j\nth}
	\leq 
	\varepsilon^{n y_{\text{min}}\nth
	},
	$$
where $y_j\nth\geq y_{\text{min}}\nth=\min_{i=1,\dots,d}y_i\nth$. 
Using the inequality above, the fact that the integral of $\tilde{p}$ is bounded by $1$,  and Stirling's formula yield
	\begingroup
    \allowdisplaybreaks
	\begin{align*}
	\int_{\domain\setminus\domain^\varepsilon}
	f_{D\nth}(\x)\tilde{p}(\x)d\x
	&\leq
	\frac{\varepsilon^{n y_{\text{min}}\nth}}{B\left(\y\nth+\boldsymbol{1}\right)}
	\int_{\domain\setminus\domain^\varepsilon} \tilde{p}(\x)d\x
	\\
	&\leq
	\frac{\varepsilon^{n y_{\text{min}}\nth} \ 
	\Gamma\left(n\norm{\y}\nth+d\right)}{\prod_{i=1}^d \Gamma\left(ny_i\nth+1\right)}
	\\
	&\sim
	(2\pi)^{-\frac{d-1}{2}} \ \varepsilon^{n y_{\text{min}}\nth}\ 
	\frac{\left(n\norm{\y\nth}+d\right)^{n\norm{\y\nth}+d-\frac{1}{2}}
	}{
	\prod_{i=1}^{d}\left(ny_i\nth+1\right)^{ny_i\nth+\frac{1}{2}}
	}
	\\
	&\sim
	\left[
	\varepsilon^{y_{\text{min}}} \prod_{i=1}^{d}\left(\frac{\norm{\y}}{y_i}\right)^{y_i}
	\right]^n
	n^{\frac{d-1}{2}}
	\left(\frac{\norm{\y}}{2\pi}\right)^{\frac{d-1}{2}}
	\left(\prod_{i=1}^{d}\frac{\norm{\y}}{y_i}\right)^{\frac{1}{2}},
	\end{align*}
	\endgroup
where $y_{\text{min}}=\min_{i=1,\dots,d}y_i>0 $.
The expression in the last display converges to zero by choosing $\varepsilon<\prod_{i=1}^{d}\left(\frac{y_i}{\norm{\y}}\right)^{\frac{y_i}{y_{\text{min}}}} $, and thus the second integral in \eqref{eq:intsplit} converges to zero, as $n\to\infty$, completing the proof.
\end{proof}
\end{theorem}

By applying Theorem \ref{thm:kconv}, it is 
straightforward to show that the asymptotic decay of the sampling probabilities is indeed polynomial as expected.  
\begin{theorem}
\label{thm:asymptoticp}
Let $p$  be the sampling probability \eqref{eq:multinomialdraw} of the block counting process of the ancestral selection graph. 
Let $\tilde{p}$ be the stationary density of the Wright-Fisher diffusion \eqref{eq:sdeWF}, assuming the 
mutation probability matrix $P$ is irreducible.
Let $\y\nth\in \frac{1}{n}\mathbb{N}^d\setminus\{\boldsymbol{0}\}$, such that $\y\nth\to \y\in\mathbb{R}^d_{>0}$, as $n\to\infty$. Then
	$$
	\lim_{n\to\infty}
	n^{d-1}\  p(n\y\nth)  = 
	\norm{\y}^{1-d} \ 
	\tilde{p}\left(\frac{\y}{\norm{\y}}\right).
	$$
\begin{proof}
Since
	$
	p(n\y\nth)=\binom{ n\norm{\y\nth}}{ n\y\nth} k(n\y\nth)
	$,
rewrite
	\begin{align*}
	n^{d-1} p(n\y\nth) =
	n^{d-1} \binom{ n\norm{\y\nth}}{ n\y\nth}B\left(n\y\nth+\boldsymbol{1}\right) \frac{k\left(n\y\nth\right)}{B\left(n\y\nth+\boldsymbol{1}\right)}.
	\end{align*}
Then, note that 
	$$
	n^{d-1} \binom{ n\norm{\y\nth}}{ n\y\nth}
	B\left(n\y\nth+\boldsymbol{1}\right)
	=
	n^{d-1} 
	\frac{\Gamma\left(n\norm{\y\nth}+1\right)}{\Gamma\left(n\norm{\y\nth}+d\right)}
	\to
	\norm{\y}^{1-d},
	$$
as $n \to \infty$, 	
whereas, by Theorem \ref{thm:kconv},     
    $
    \frac{k\left(n\y\nth\right)}{B\left(n\y\nth+\boldsymbol{1}\right)}
    $ 
converges to 		
	$
	\tilde{p}\left(\frac{\y}{\norm{\y}}\right)
	$. This completes the proof.
\end{proof}
\end{theorem}

%%%%%%%%%%%%%%%%%%%%%%%%%%%%%%%%%%%%%%%%%%%%%%%%%%%%%%%%%%%%%%%%%%%%%%%%%%%%%%%%%%%

\section{Asymptotic behaviour of the transition probabilities}
\label{sect:probabilities}
The goal of this section is to study the asymptotic behaviour of $p(n\y\nth-\mathbf{v}\mid n\y\nth)$, the transition probabilities \eqref{defrho} of the block counting jump chain $\kc$ of the ASG, as $n\to\infty, \y\nth\to\y$.
%In order to obtain a significant limit, the following scaling is needed. 
By letting $\kc\nth$, $n\in\mathbb{N},$ be independent copies of $\kc$, and  $\Y\nth=\frac{1}{n}\kc\nth\subset \frac{1}{n}\mathbb{N}^d\setminus\{\boldsymbol{0}\}$,
the asymptotic behaviour of the transition probabilities can be interpreted as the limit of the  transition probabilities of $\Y\nth$, $\rho\nth(\mathbf{v}\mid \y\nth)=p(n\y\nth-\mathbf{v}\mid n\y\nth)$.

As for the sampling probabilities, 
when mutations are parent dependent, an explicit expression for the backward transition probabilities is not available, in fact  expression \eqref{defrho} is written in terms of the unknown probability $\pi$.

In the PIM neutral case, $P_{ij}=Q_j, \gamma_j=0, \  i,j=1,\dots,d$, the sampling probabilities can be explicitly written as
    $
    \pi[i|\n]=  
    \frac{ n_i+\theta Q_i}{\norm{\n}+\theta}
    $, 
see e.g. \cite{stephens2000}.
In this case thus  $ \pi[i|n\y\nth] $  converges to $\frac{y_i}{\norm{\y}}$. It turns out that in the general case the limit of the probabilities $\pi$ is the same. 
To prove this, $\pi$ is written in 
terms of  the function $k$, defined in \eqref{eq:kdef},  by using \eqref{defpi} and \eqref{eq:p_intermsof_k}, 
    \begin{align*}
    \pi[i|\n]=  
    \frac{k(\n+\e_i)}{k(\n)}.
    \end{align*}
The results of the previous section, combined with the expression above, make it straightforward to study the asymptotic behaviour of $\pi[i|n\y\nth]$, and, consequently, of the transition probabilities. 
\begin{proposition}
\label{prop:asymptoticpi}
Let $\pi$ be defined as in \eqref{defpi}. Assume the mutation probability matrix $P$ is irreducible.
Let $\y\nth\in \frac{1}{n}\mathbb{N}^d\setminus\{\boldsymbol{0}\}$, such that $\y\nth\to \y\in\mathbb{R}^d_{>0}$, as $n\to\infty$. Then, for $i=1,\dots,d, $
	\begin{align*}
	\lim_{n\to\infty} \pi[i|n\y\nth]= \frac{y_i}{\norm{\y}}.
	\end{align*}
\begin{proof}
Rewrite
	$$
	\pi[i|n\y\nth]=
	\frac{k(n\y\nth+\e_i)}
	{B(n\y\nth+\e_i+\boldsymbol1)}
	\frac{B(n\y\nth+\boldsymbol1)}
	{k(n\y\nth)}
	\frac{B(n\y\nth+\e_i+\boldsymbol1)}
	{B(n\y\nth+\boldsymbol1)}.
	$$
By Theorem \ref{thm:kconv}, 
	$
	\frac{k(n\y\nth+\e_i)}
	{B(n\y\nth+\e_i+\boldsymbol1)}
	$ 
and 
	$
	\frac{k(n\y\nth)}
	{B(n\y\nth+\boldsymbol1)}
	$
both converge to $\tilde{p}\left(\frac{\y}{\norm{\y}}\right)$, as $n\to\infty$. Calculating,
	$$
	\frac{B(n\y\nth+\e_i+\boldsymbol1)}
	{B(n\y\nth+\boldsymbol1)}
	=
	\frac{\Gamma(n y_i\nth +2)}{\Gamma(n y_i\nth +1)}
	\frac{\Gamma(n\norm{\y\nth}+d)}{\Gamma(n\norm{\y\nth}+1+d)}
	= \frac{n y_i\nth +1}{n \norm{\y\nth}+d},
	$$
and letting $n \to \infty$ concludes the proof.
\end{proof}
\end{proposition}
Finally, knowing the asymptotic behaviour of $\pi$  directly solves the problem of analysing the asymptotic behaviour of the  transition probabilities.

\begin{corollary}
\label{coroll:asymptoticrho}
Let $\rho\nth(\mathbf{v}\mid\y\nth)=p(n\y\nth-\mathbf{v}\mid \y\nth)$ be the transition probabilities defined in \eqref{defrho}. Under the assumptions of Proposition \ref{prop:asymptoticpi},
for $i,j=1,\dots,d, $
	\begin{align*}
	&\lim_{n\to\infty}
	\rho\nth(\e_j\mid\y\nth)= \frac{y_j}{\norm{\y}} ;
	\\
	&
	\lim_{n\to\infty}     n\rho\nth(\e_j-\e_i\mid\y\nth)=\frac{\theta P_{ij}y_i}{\norm{\y}^2} ;
	\\
	&
	\lim_{n\to\infty}     n\rho\nth(-\e_j\mid\y\nth)=\frac{ |\gamma_j|y_j}{\norm{\y}^2}. 
	\end{align*}
\begin{proof}
By Proposition \ref{prop:asymptoticpi},
$\pi[i|n\y\nth-\e_j]\to \frac{y_i}{\norm{\y}}$,
$\pi[j|n\y\nth]\to \frac{y_j}{\norm{\y}}$,
for all $i,j=1,\dots,d$.
Using these limits, together with basic limit calculations, in \eqref{defrho}, gives the result. 
\end{proof}
\end{corollary}

%%%%%%%%%%%%%%%%%%%%%%%%%%%%%%%%%%%%%%%%%%%%%%%%%%%%%%%%%%%%%%%%%%%%%%%%%%%%%%%%%%%

\section*{Appendix: central and local limit theorems for a sequence of Dirichlet random vectors}

This section is devoted to the study of the asymptotic behaviour of  the sequence of Dirichlet distributed random vectors presented in Proposition \ref{prop:convdirich}. 
While this section can be considered as a standard, yet not obvious, exercise in probability theory, 
it is essential to derive central and local limit theorems for the specific sequence of random vectors in this paper.

Let $\mathbf{D}\nth=(D_1\nth,\dots,D_d\nth)\in \mathcal{S}$  be a Dirichlet distributed random vector with concentration parameters
$\boldsymbol{\alpha}\nth\in \mathbb{R}_{>0}^d$ such that 
    \begin{equation}
    \label{eq:alphalinear}
    \lim_{n\to\infty}\frac{\boldsymbol{\alpha}\nth}{n}=\boldsymbol{\alpha}\in \mathbb{R}_{>0}^d.
    \end{equation}

A central limit theorem for the sequence $\mathbf{D}\nth$ is obtained using the well-known,
see e.g. \cite{Devroye1986}%Thm 4.1,Ch XI
, relationship between Dirichlet and Gamma distributions,
    \begin{align}
    \label{gamma2dirichlet}
    \mathbf{D}\nth \stackrel{d}{=} \frac{\mathbf{G}\nth}{\norm{\mathbf{G}\nth}},    
    \end{align}
where $\mathbf{G}\nth=(G_1\nth,\dots,G_d\nth)$ is a vector of independent Gamma  random variables with shape parameters $\alpha_i\nth, i=1,\dots,d$ and rate parameter $\beta\in \mathbb{R}_{>0}$.
Note that  $\beta$ is irrelevant in the transformation from $\mathbf{G}\nth$ to $\mathbf{D}\nth$.
Furthermore, since the Gamma distribution is infinitely divisible, it is possible to write each of the Gamma random variables as a sum of independent Gamma random variables,
and thus 
for each variable $G_i\nth$ the following central limit theorem  holds
    $$
    \frac{\sqrt{n}\beta}{\sqrt{\alpha_i\nth}}\left(
    \frac{1}{n}G_i\nth -  \frac{\alpha_i\nth}{n\beta}
    \right) \xrightarrow[n\to\infty]{d}
    \mathcal{N}
    \left(0,1 \right),
    \quad i=1,\dots,d.
    $$
Consequently, by \eqref{eq:alphalinear} and  
independence of the components of the random vectors $\mathbf{G}\nth$, the following central limit theorem for $\mathbf{G}\nth$ holds	
    $$
    \sqrt{n}\left(
    \frac{1}{n}\mathbf{G}\nth -  \frac{\boldsymbol{\alpha}}{\beta}
    \right) \xrightarrow[n\to\infty]{d}
    \mathcal{N}_d
    \left(\boldsymbol{0},\frac{1}{\beta^2}diag(\boldsymbol{\alpha})\right),
    $$
where $diag(\boldsymbol{\alpha})$ is the diagonal matrix with $\boldsymbol{\alpha}$ on the diagonal.
Since 
$\mathbf{D}\nth$ can be written as  a function of $\mathbf{G}\nth $,  
%$\mathbf{D}\nth \stackrel{d}{=} f(\mathbf{G}\nth)$, as in 
as in \eqref{gamma2dirichlet},
applying the multivariate delta method 
%( Sokol Ronn 3.7.6) 
yields
    $$
    \sqrt{n}\left(
    \mathbf{D}\nth 
    -  \frac{\boldsymbol{\alpha}}{\norm{\boldsymbol{\alpha}}}
    \right) \xrightarrow[n\to\infty]{d}
    \mathcal{N}_d
    \left(\boldsymbol{0},\Sigma(\boldsymbol{\alpha})\right),
    $$
with 
    $$
    \Sigma(\boldsymbol{\alpha})=
    J \left(\frac{\boldsymbol{\alpha}}{\beta}\right)
    \frac{1}{\beta^2}\text{diag}(\boldsymbol{\alpha})
    J \left(\frac{\boldsymbol{\alpha}}{\beta}\right)^T, 
    $$
where $J$ is the Jacobian matrix associated to the transformation in \eqref{gamma2dirichlet}. 
Calculating the Jacobian and  multiplying the matrices, which is omitted,
yields
    $$
    \Sigma_{ij}(\boldsymbol{\alpha})= 
    \frac{\alpha_i}{\norm{\boldsymbol{\alpha}}^3}
    (\delta_{ij}\norm{\boldsymbol{\alpha}}-\alpha_j),
    \quad i,j=1,\dots,d .
    $$
As expected, the covariance matrix above does not depend on the auxiliary rate parameter $\beta$. 

Note that, if $\norm{\boldsymbol{\alpha}}=1$, then
$\Sigma(\boldsymbol{\alpha})=\sigma(\boldsymbol{\alpha})$, where $\sigma$ is the diffusion matrix in \eqref{eq:sdeWF}.
Consequently, the Wright-Fisher diffusion matrix can be interpreted as the covariance matrix of the Gaussian limit of a sequence of Dirichlet random vectors. 

The Gaussian limiting vector,
being the limit of a sequence of Dirichlet random vectors,  has a degenerate distribution, as its last component can be expressed in terms of the first $d-1$ components. Therefore, 
in order to work with density functions and obtain a local limit result, the state space is restricted to $\mathbb{R}^{d-1}$ by excluding the last component of each vector in the remaining part of this section. 

Let $\phi$, be the pdf of the $d-1$ dimensional centred Gaussian vector with covariance matrix $\Sigma_{d-1}(\boldsymbol{\alpha})$, the restriction of $\Sigma(\boldsymbol{\alpha})$ to the first $d-1$ components. 
Let 
    $$
    f_{\mathbf{D}\nth}(\x)
    =\frac{1}{B(\boldsymbol{\alpha}\nth)}
    \prod_{i=1}^d x_i^{\alpha_i\nth-1},
    \quad 
    \mathbf{x}\in\domain,
    $$
be the pdf of (the first $d-1$ components)  of  $\mathbf{D}\nth$,
where $x_d$ stands for $1-\sum_{i=1}^{d-1}x_i$.
Then the pdf of 
    $
    \sqrt{n}\left(
    \mathbf{D}\nth 
    -  \frac{\boldsymbol{\alpha}\nth}{\norm{\boldsymbol{\alpha}\nth}}
    \right)
    $
is  given by
    $$
    \phi_n (\mathbf{u}) =
    n^{-\frac{1}{2}(d-1)}f_{\mathbf{D}\nth}
    \left(\frac{1}{\sqrt{n}}\mathbf{u} 
    +\frac{\boldsymbol{\alpha}\nth}{\norm{\boldsymbol{\alpha}\nth}}\right),
    \quad \text{ for  }\mathbf{u}\in -\sqrt{n}
    \left(\domain
    -  \frac{\boldsymbol{\alpha}\nth}{\norm{\boldsymbol{\alpha}\nth}}
    \right), 
    $$
and equal to $0$ otherwise. In general, convergence in distribution does not imply convergence of probability density functions,
however, 
under the conditions of the converse to Scheffe's theorem, see \cite{boos1985}, it does.
More precisely, if the sequence $\phi_n$ is bounded and uniformly equicontinuous,  then  
$
\phi_n \to \phi
$
, as $n \to\infty$, uniformly on $\mathbb{R}^{d-1}$. In the remaining part of this section, the conditions for the convergence of densities are verified, again, it is fundamental that the parameters of the Dirichlet vectors grow to infinity linearly,  as assumed in \eqref{eq:alphalinear}. 

In order to show boundedness of $\phi_n$, 
first notice that, for a fixed $n$, the maximum of $\phi_n$ is reached at 
    $
    %\mathbf{u}=
    \sqrt{n}\left(
    \frac{\boldsymbol{\alpha}\nth-\boldsymbol{1}}{\norm{\boldsymbol{\alpha}\nth}-d} -
    \frac{\boldsymbol{\alpha}\nth}{\norm{\boldsymbol{\alpha}\nth}}
    \right),
    $
and thus 
    $$
    ||\phi_n||_\infty=
    \sup_{\mathbf{u}\in \mathbb{R}^{d-1}}
    \phi_n(\mathbf{u})=
    n^{-\frac{d-1}{2}}
    \frac{1}{B(\boldsymbol{\alpha}\nth)}
    \prod_{i=1}^d 
    \left( 
    \frac{\alpha_i\nth-1}{\norm{\boldsymbol{\alpha}\nth}-d}
    \right)^{\alpha_i\nth-1}.
    $$
Using  Stirling's formula for the Beta function, 
    $$
    ||\phi_n||_\infty \sim
    \left(
    (2\pi)^{{d-1}}\prod_{i=1}^{d}\frac{\alpha_i\nth}{\norm{\boldsymbol{\alpha}\nth}}
    \right)^{-\frac{1}{2}}
    \left(\frac{\norm{\boldsymbol{\alpha}\nth}}{n}
    \right)^{\frac{1}{2}(d-1)}.
    $$
From the expression above, it is clear that 
if $\norm{\boldsymbol{\alpha}\nth}$ grows faster than linearly, the sequence is not bounded and if it grows slower it converges to zero. By assumption \eqref{eq:alphalinear},
$\norm{\boldsymbol{\alpha}\nth}$ grows linearly and therefore
$||\phi_n||_\infty$ converges to
    \begin{equation}
    \label{eq:norm_of_phi}
    \left(
    (2\pi)^{d-1}
    \prod_{i=1}^{d}\frac{\alpha_i}{\norm{\boldsymbol{\alpha}}}
    \frac{1}{\norm{\boldsymbol{\alpha}}^{d-1}}
    \right)^{-\frac{1}{2}},
     \end{equation}
as  $n\to\infty$, 
and thus the sequence of density functions $\phi_n$ is bounded.

Following the same type of argument as for the sequence of density functions, it is straightforward to show that also each sequence of their first order partial derivatives is bounded. Furthermore, the density functions $\phi_n$ are smooth on a compact set. Therefore, the sequence of densities is uniformly Lipschitz continuous,
% i.e.
%$\forall \mathbf{u},\mathbf{u}' \in \mathbb{R}^{d-1}, \quad |\phi_n(\mathbf{u})-\phi_n(\mathbf{u}')|\leq L |\mathbf{u}-\mathbf{u}'|$,
which verifies the uniform equicontinuity condition, proving the local limit result as stated in  Proposition \ref{prop:convdirich}. 

This section is concluded with the following observation. 

\begin{Remark}
The argument in this section  provides a straightforward derivation of an expression for the determinant of the diffusion matrix $\sigma_{d-1}(\x)$ of the Wright-Fisher diffusion in \eqref{eq:sdeWF}, as explained in the following.
We have proven that  $||\phi_n-\phi||_\infty\to 0$,  as $n\to \infty$,
which implies $||\phi_n ||_\infty \to ||\phi||_\infty$. Note that $
    ||\phi||_\infty
    =\left[(2 \pi)^{d-1} 
    \det(\Sigma_{d-1}(\boldsymbol{\alpha}))
    \right]^{-1/2}$.
Furthermore, we have proven that  $||\phi_n||_\infty$ converges to \eqref{eq:norm_of_phi}.
Therefore, it must be that  $||\phi||_\infty$ is equal to \eqref{eq:norm_of_phi}, which  implies
  $$
    \det(\Sigma_{d-1}(\boldsymbol{\alpha}))=
    \prod_{i=1}^{d}\frac{\alpha_i}{\norm{\boldsymbol{\alpha}}}
    \frac{1}{\norm{\boldsymbol{\alpha}}^{d-1}}.
    $$
 Therefore, the determinant of the diffusion matrix of the Wright-Fisher diffusion, for $\x\in\mathcal{S}$, is
    $$
    \det(\sigma_{d-1}(\x))=
    \prod_{i=1}^{d}x_i.
    $$   
 \end{Remark}

%%%%%%%%%%%%%%%%%%%%%%%%%%%%%%%%%%%%%%%%%%%%%%%%%%%%%%%%%%%%%%%%%%%
%%                                                               %%
%% Supplementary Material, if any, should be provided in         %%
%% {supplement} environment  with title and short description.   %%
%%                                                               %%
%%%%%%%%%%%%%%%%%%%%%%%%%%%%%%%%%%%%%%%%%%%%%%%%%%%%%%%%%%%%%%%%%%%

%\begin{supplement}
%\stitle{Title of Supplement A.}
%\sdescription{Short description of Supplement A.}
%\end{supplement}
%\begin{supplement}
%\stitle{Title of Supplement B.}
%\sdescription{Short description of Supplement B.}
%\end{supplement}

%%%%%%%%%%%%%%%%%%%%%%%%%%%%%%%%%%%%%%%%%%%%%%%%%%%%%%%%%%%%%%%%%%%
%%                                                               %%
%% Use the two commands below for producing your bibliography    %%
%% with bibtex, then comment again the commands and include the  %%
%% content of the .bbl file in this file below the commands.     %%
%%                                                               %%
%%%%%%%%%%%%%%%%%%%%%%%%%%%%%%%%%%%%%%%%%%%%%%%%%%%%%%%%%%%%%%%%%%%

%\bibliographystyle{amsplain}
%\bibliography{yourbibfilename}

% add below the content of your .bbl file produced by bibtex.

%%%%%%%%%%%%%%%%%%%%%%%%%%%%%%%%%%%%%%%%%%%%%%%%%%%%%%%%%%%%%%%%%%%
%%                                                               %%
%% You may add acknowledgments (optional).                       %%
%%                                                               %%
%%%%%%%%%%%%%%%%%%%%%%%%%%%%%%%%%%%%%%%%%%%%%%%%%%%%%%%%%%%%%%%%%%%

\ACKNO{The authors would like to thank the anonymous reviewer whose comments led to an improvement of the manuscript.}

%%%%%%%%%%%%%%%%%%%%%%%%%%%%%%%%%%%%%%%%%%%%%%%%%%%%%%%%%%%%%%%%%%%
%%                                                               %%
%% You have reached the end of your document.                    %%
%%                                                               %%
%%%%%%%%%%%%%%%%%%%%%%%%%%%%%%%%%%%%%%%%%%%%%%%%%%%%%%%%%%%%%%%%%%%

\end{document}